\documentclass[12pt]{amsart}
\usepackage{natbib}
 \bibpunct[, ]{(}{)}{,}{a}{}{,}

\begin{document}

\title{Solving maximum cut problems by simulated annealing}
\author{Tor Myklebust}
\thanks{Tor Myklebust:
Department of Combinatorics and Optimization, Faculty
of Mathematics, University of Waterloo, Waterloo, Ontario N2L 3G1,
Canada (e-mail: tmyklebu@csclub.uwaterloo.ca).  Research of this author
was supported in part by an NSERC Doctoral Scholarship.}
\date{May 8, 2015}

\begin{abstract}
  This paper gives a straightforward implementation of simulated annealing for
  solving maximum cut problems and compares its performance to that of some
  existing heuristic solvers.  The formulation used is classical, dating to
  a 1989 paper of Johnson, Aragon, McGeoch, and Schevon.
  This implementation uses no structure peculiar to the maximum cut problem,
  but its low per-iteration cost allows it to find better solutions than were
  previously known for 40 of the 89 standard maximum cut instances tested
  within a few minutes of computation.
\end{abstract}

\maketitle

\begin{section}{Introduction}
  Simulated annealing is a classical heuristic method for solving optimisation
  problems due independently to Kirkpatrick, Gelatt, and Vecchi
  \citep{KirkpatrickGelattVecchi} and \v{C}ern\'{y} \citep{Cerny}.  The maximum
  cut problem is a famous NP-hard optimisation problem.  The purpose of the
  present work is to compare a straightforward implementation of simulated
  annealing with the computational state of the art for heuristic solution to
  the maximum cut problem as reported by Mart\'{i}, Duarte, and Laguna
  \citep{MartiDuarteLaguna}.  This work also compares with Burer, Monteiro, and
  Zhang's CirCut heuristic \citep{BurerMonteiroZhang} and Festa, Pardalos,
  Resende, and Ribeiro's Variable Neighbourhood Search with Path Relinking
  heuristic \citep{FestaPardalosResendeRibeiro}.

  The maximum cut problem is defined on a graph $G = (V,E)$ with weights
  $w : E \to \mathbb{Z}$ on its edges.  It asks for the vertex set
  $S \subseteq V$ maximising
  $\sum_{uv : u \in S, v \in V \setminus S, uv \in E} w(uv)$, the total weight
  of edges with exactly one end in $S$.  $S$ and $V \setminus S$ here are
  called the shores of the cut.
  
  The maximum cut problem is NP-hard \citep{Karp}---under the assumption
  that $P \neq NP$, no polynomial-time algorithm exists for computing an
  optimal solution to the maximum cut problem.  More is true; H\aa stad proved
  \citep{Hastad} that, under the assumption that $P \neq NP$, no
  polynomial-time algorithm exists that always computes a solution whose
  objective value is within a factor $0.942$ of the optimum.

  Goemans and Williamson \citep{GoemansWilliamson} devised a polynomial-time
  randomised rounding algorithm based on an SDP relaxation that computes, with
  nonzero, constant probability, a solution whose objective value is within a
  factor $\alpha \simeq 0.87856$ of the optimum whenever $w$ only takes
  positive values.  Khot, Kindler, Mossel, and O'Donnell \citep{Khot} proved
  that, if the Unique Games Conjecture is true, then even in the case where $w$
  only takes positive values there is no polynomial-time algorithm that always
  computes a solution whose objective value is within a constant factor larger
  than $\alpha$ of the optimum.

  Johnson, Aragon, McGeoch, and Schevon \citep{JohnsonAragon1} gave an
  empirical evaluation of simulated annealing applied to a certain graph
  partitioning problem---find the maximum cut in a graph with both shores
  having the same number of vertices.

  Burer, Monteiro, and Zhang \citep{BurerMonteiroZhang} proposed approximately
  solving a certain nonconvex continuous relaxation of the maximum cut problem
  and then using randomised rounding on the solution found.  With some
  refinements, this approach is implemented in the CirCut code of Y. Zhang.
\end{section}

\begin{section}{Simulated annealing}\label{sec:simann}
  Simulated annealing works by taking a random walk on the space of feasible
  solutions to a problem.  The random walk considers a random local
  modification of the current solution.  If that modification leads to a better
  solution, it is always made.  If not, it is made with some probability
  that decreases exponentially with the amount of objective degradation.
  Typically, the parameter of this exponential distribution is increased as the
  number of iterations increases (thus reducing the probability of accepting a
  step that worsens the objective value).  This parameter, as a function of the
  iteration count, is called the ``annealing schedule.''

  One can represent the space of feasible solutions to an instance the maximum
  cut problem defined on the graph $G = (V,E)$ as the space $\{-1, 1\}^V$ of
  functions that specify which shore of the cut each vertex lies on.  Here, two
  cuts are adjacent---a single step of the random walk can get from one to the
  other---if they differ in the placement of at most one vertex.  This setup
  allows an especially fast implementation of simulated annealing, since the
  objective function change that would be caused by moving a vertex $v$ to
  the other shore changes only when $v$ or one of its neighbours is moved.
  (This is the same setup used by \citet{JohnsonAragon1} in their study of
  simulated annealing for the closely-related graph partitioning problem.)

  The C++ code in Figure \ref{code} is an implementation of simulated annealing
  for the maximum cut problem for graphs with integer edge weights and at most
  33333 vertices.  It accepts input in the standard format---one line telling
  the number of vertices and the number of edges of the graph, followed by one
  line for each edge giving the edge's two endpoints and its weight.
  
  The \texttt{heatmax} parameter and the increment \texttt{2e-6} were selected
  because they gave acceptable results in an amount of time broadly comparable
  to that taken by Mart\'{i} et al.'s Scatter Search heuristic on larger
  graphs.  (Whereas Mart\'{i} et al. take up to 23 minutes on their hardware on
  some instances, all runs of simulated annealing reported here took
  considerably less than ten minutes.) Since this annealing schedule is
  oblivious to the size and sparsity of the input graph and the magnitudes of
  the weights on its edges, it is quite likely that a more effective
  general-purpose annealing schedule can be found.

  Note that each iteration of this simulated annealing heuristic takes either a
  very small constant amount of time (two random numbers, one exponential, and
  a little arithmetic) or time linear in the degree of the selected vertex.
  This allows the code to explore a substantially larger number of feasible
  solutions---about $5 \times 10^9$ here---than other approaches might within
  the same timeframe.

\begin{figure}
\begin{verbatim}
#define FOR(i,n) for (int i=0;i<(n);i++)
vector<pair<int,int> > edges[33333];
int n, m, side[33333], chg[33333];

int main() {
  if (2 != scanf("%i%i", &n, &m)) abort();
  while(m--) {
    int a, b, w;
    if (3 != scanf("%i%i%i", &a, &b, &w)) abort();
    a--;b--;
    edges[a].push_back(make_pair(b, w));
    edges[b].push_back(make_pair(a, w));
    chg[a] += w; chg[b] += w;
  }
  FOR(i,n) side[i] = 1;
  int obj = 0, best = 0;
  double heatmax = 10000;
  for (double heat = 0; heat < heatmax; heat += 2e-6) {
    int k = rand() % n;
    if (!(drand48() > exp(heat*chg[k]/best))) {
      obj += chg[k]; side[k] = -side[k]; chg[k] = -chg[k];
      for (auto e : edges[k]) {
        int v = e.first;
        chg[v] += e.second * (2 - 4 * (side[k] != side[v]));
      }
      if (obj > best) {
        FOR(i,n) printf(" %i", side[i]); printf("\n");
        best = obj;
      }
    }
  }
}
\end{verbatim}
\caption{A simulated annealing heuristic for maximum cut in C++.}
\label{code}
\end{figure}
\end{section}

\begin{section}{Computational experiments}
  The simulated annealing code of Figure \ref{code} was run on a computer with
  an Intel Xeon E3-1245v2 processor (8 logical cores) and 32 GB of RAM.  This
  computer runs Linux.  The code of Figure \ref{code} was compiled with
  optimisation (\texttt{-O3}) using GCC 4.9.2.
  
  CirCut 1.0612 \citep{BurerMonteiroZhang} was run on the same hardware, but
  compiled according to the Makefile distributed with the software.  I changed
  the parameters \texttt{npert} to 300 and \texttt{multi} to 25 from the default
  in order to get better solutions on average.  This comes at a substantially
  greater cost in running time.

  I was able to build the maximum cut heuristics, including VNSPR, available
  from M. Resende's website\footnote{At the time of writing,
  \texttt{http://mauricio.resende.info/src/maxcut-heuristics.tar.gz}.}, but I
  could not get VNSPR to produce good results in ten minutes of computation
  time.  Solution values and times for VNSPR reported in this paper are taken
  from the tables linked on Optsicom's website\footnote{At the time of writing,
  \texttt{http://www.optsicom.es/maxcut}.}, and these very closely resemble
  those reported by Festa et al. in \citep{FestaPardalosResendeRibeiro}.  The
  computation times for VNSPR listed here should be assumed to have been
  recorded on a far slower machine---an SGI Challenge from the 20th
  century---than the times given for other heuristics.
  
  I requested from R. Mart\'{i} the implementations of Scatter Search (SS)
  \citep{MartiDuarteLaguna}, CirCut \citep{BurerMonteiroZhang}, and Variable
  Neighbourhood Search with Path Relinking (VNSPR)
  \citep{FestaPardalosResendeRibeiro} used in the experiments in
  \citep{MartiDuarteLaguna}, but my request went unfulfilled.  However, Y. Zhang
  has made available via his website the source code for CirCut and M. Resende
  has made available via his website the source code for all of the heuristics
  evaluated by Festa et al.  The solution values and computation times reported
  for Scatter Search and VNSPR are taken from the tables linked on Optsicom's
  website accompanying the paper by Mart\'{i}, Duarte, and Laguna.

  Mart\'{i} et al. report that the hardware they used for Scatter Search is
  roughly comparable to (but slightly weaker than) that used for this simulated
  annealing code.

  The comparison was run on three sets of test cases, all retrievable via
  Optsicom's website.  Following Mart\'{i}, Duarte, and Laguna, these are
  called Sets 1, 2, and 3.  Set 1 was generated by Helmberg and Rendl
  \citep{HelmbergRendl} using G. Rinaldi's graph generator \texttt{rudy}.  Set 2
  was generated by Festa, Pardalos, Resende, and Ribeiro
  \citep{FestaPardalosResendeRibeiro} from Ising spin glass models.  Set 3 has
  four graphs from the 7th DIMACS Implementation Challenge. 

  Results are given in Figures \ref{set1}, \ref{set2}, and \ref{set3}.  The
  first column of each row gives the name of the graph.  The columns labelled
  SS, CirCut, VNSPR, and SA give the objective value of the best solution found
  by each heuristic.  The adjacent column labelled Time gives the computation
  time, in seconds, required by that heuristic.

  A summary of the results is given in Figure \ref{summarytable}.  The first
  column gives the number of the dataset in question (1, 2, or 3).  The
  columns labelled SS, CirCut, VNSPR, and SA give the number of graphs in that
  dataset on which each heuristic found the best solution of the four.
  The adjacent column labelled ``only'' gives the number of graphs in that
  dataset on which each heuristic found strictly better solutions than the
  other three heuristics.

  \begin{figure}
  \begin{scriptsize}
  \begin{tabular}{l|rr|rr|rr|rr}
    Graph & SS & Time & CirCut & Time & VNSPR & Time & SA & Time\\
    \hline
g1&11624&139.0&11624&243.378&11621&22732.0&11621&294.680\\
g2&11620&167.2&11600&230.651&11615&22719.0&11612&327.224\\
g3&11622&180.1&11621&220.899&11622&23890.0&11618&295.132\\
g4&11646&194.4&11633&186.158&11600&24050.0&11644&294.280\\
g5&11631&205.2&11630&322.106&11598&23134.0&11628&299.572\\
g6&2165&175.7&2178&209.259&2102&18215.4&2178&246.516\\
g7&1982&175.5&2006&251.745&1906&17715.8&2006&204.756\\
g8&1986&194.9&2004&204.970&1908&19333.8&2005&206.232\\
g9&2040&158.0&2052&223.100&1998&15224.5&2054&205.500\\
g10&1993&209.7&1996&217.800&1910&16268.9&1999&204.992\\
g11&562&171.8&558&63.512&564&10084.0&564&188.928\\
g12&552&241.5&554&38.496&556&10852.0&554&188.712\\
g13&578&227.5&578&103.056&580&10749.0&580&194.996\\
g14&3060&186.5&3057&97.148&3055&16734.0&3063&252.456\\
g15&3049&142.8&3037&73.096&3043&17184.0&3049&220.432\\
g16&3045&161.9&3047&86.183&3043&16562.0&3050&219.472\\
g17&3043&312.8&3043&125.603&3030&18554.6&3045&219.464\\
g18&988&174.3&988&79.426&916&12578.3&990&234.828\\
g19&903&128.3&902&71.282&836&14545.6&904&196.020\\
g20&941&191.0&940&77.651&900&13325.6&941&194.644\\
g21&930&233.1&929&93.120&902&12884.6&927&194.656\\
g22&13346&1335.8&13352&457.417&13295&197654.0&13158&294.572\\
g23&13317&1021.7&13324&363.187&13290&193707.0&13116&287.836\\
g24&13303&1191.0&13319&388.584&13276&195749.0&13125&289.108\\
g25&13320&1299.2&13333&324.094&12298&212563.3&13119&315.668\\
g26&13294&1415.1&13313&467.400&12290&228969.4&13098&288.544\\
g27&3318&1437.5&3326&368.070&3296&35652.3&3341&214.312\\
g28&3285&1314.0&3286&344.389&3220&38654.8&3298&251.784\\
g29&3389&1265.8&3390&411.182&3303&33694.5&3394&214.440\\
g30&3403&1196.3&3405&479.210&3320&34457.5&3412&214.628\\
g31&3288&1336.2&3293&372.585&3202&36658.0&3309&214.044\\
g32&1398&900.6&1390&153.000&1396&82345.0&1410&193.964\\
g33&1362&925.6&1354&119.342&1376&76282.0&1376&194.092\\
g34&1364&950.2&1370&194.402&1372&79406.0&1382&194.116\\
g35&7668&1257.5&7672&335.542&7635&167221.0&7485&262.772\\
g36&7660&1391.9&7656&264.945&7632&167203.0&7473&265.460\\
g37&7664&1386.8&7672&340.885&7643&170786.0&7484&287.644\\
g38&7681&1011.5&7674&364.396&7602&178569.8&7479&263.840\\
g39&2393&1310.9&2393&264.634&2303&42584.2&2405&209.412\\
g40&2374&1166.4&2387&279.101&2302&39548.9&2378&208.376\\
g41&2386&1016.5&2398&239.882&2298&40025.0&2405&208.136\\
g42&2457&1458.2&2465&278.497&2390&41254.5&2465&210.452\\
g43&6656&405.8&6658&167.068&6659&35324.0&6658&245.056\\
g44&6648&355.9&6649&124.051&6642&34519.0&6646&241.360\\
g45&6642&354.3&6649&127.257&6646&34179.0&6652&241.496\\
g46&6634&498.1&6634&118.105&6630&38854.0&6647&244.896\\
g47&6649&359.1&6655&106.296&6640&36587.2&6652&241.848\\
g48&6000&20.1&6000&54.245&6000&64713.0&6000&210.248\\
g49&6000&35.1&5940&85.079&6000&64749.0&6000&210.180\\
g50&5880&26.8&5876&156.548&5880&147132.0&5858&210.864\\
g51&3846&513.0&3841&99.264&3808&89965.5&3841&233.704\\
g52&3849&550.8&3840&146.937&3816&95984.5&3845&228.092\\
g53&3846&423.8&3840&103.453&3802&92458.8&3845&230.052\\
g54&3846&429.0&3850&142.984&3820&98458.0&3845&228.340\\
  \end{tabular}
  \end{scriptsize}
  \caption{Comparison on Set 1 instances.}
  \label{set1}
  \end{figure}

  \begin{figure}
  \begin{tabular}{l|rr|rr|rr|rr}
    Graph & SS & Time & CirCut & Time & VNSPR & Time & SA & Time\\
    \hline
sg3dl051000&110&1.9&110&6.178&104&1196.7&110&188.600\\
sg3dl052000&112&1.9&112&6.070&106&1196.7&112&188.132\\
sg3dl053000&106&2.1&106&7.272&102&1196.7&106&188.648\\
sg3dl054000&114&2.1&114&8.468&106&1196.8&114&188.856\\
sg3dl055000&112&2.3&112&5.318&106&1196.7&112&188.356\\
sg3dl056000&110&2.1&110&5.872&108&1196.8&110&189.532\\
sg3dl057000&112&2.0&112&4.465&108&1196.8&112&188.196\\
sg3dl058000&108&2.1&108&5.667&106&1196.8&108&188.432\\
sg3dl059000&110&1.8&110&5.544&104&1196.7&110&188.660\\
sg3dl0510000&112&1.4&112&7.105&106&1196.7&112&188.284\\
sg3dl101000&882&406.1&884&106.747&892&20409.0&890&190.900\\
sg3dl102000&894&302.4&894&67.125&900&20873.0&900&190.944\\
sg3dl103000&884&410.4&878&83.015&884&20574.0&888&191.540\\
sg3dl104000&892&485.9&894&100.869&896&19786.0&896&191.020\\
sg3dl105000&880&400.9&876&72.593&882&19160.0&884&191.100\\
sg3dl106000&870&461.8&876&119.428&880&17872.0&888&191.012\\
sg3dl107000&890&386.2&892&117.473&896&21044.0&898&191.204\\
sg3dl108000&880&466.9&876&116.292&880&19760.0&878&191.332\\
sg3dl109000&888&493.6&896&77.126&898&20930.0&898&190.880\\
sg3dl1010000&886&352.8&886&86.571&890&20028.0&894&197.388\\
sg3dl141000&2428&1320.6&2416&334.266&2416&188390.0&2442&204.244\\
sg3dl142000&2418&1121.1&2440&333.137&2416&187502.0&2454&201.704\\
sg3dl143000&2410&1215.8&2416&340.113&2406&190028.0&2438&201.728\\
sg3dl144000&2422&1237.2&2420&311.276&2418&198809.0&2444&201.828\\
sg3dl145000&2416&1122.5&2422&347.010&2416&196725.0&2438&202.732\\
sg3dl146000&2424&1213.9&2428&319.165&2420&189366.0&2440&203.872\\
sg3dl147000&2404&1230.6&2416&326.425&2404&187902.0&2432&201.508\\
sg3dl148000&2416&1132.0&2426&328.149&2418&194838.0&2436&201.556\\
sg3dl149000&2412&1213.9&2394&344.832&2384&193627.0&2420&237.256\\
sg3dl1410000&2430&1125.8&2426&315.461&2422&196454.0&2452&201.628\\
  \end{tabular}
  \caption{Comparison on Set 2 instances.}
  \label{set2}
  \end{figure}

  \begin{figure}
  \begin{tabular}{l|rr|rr}
    Graph & SS & Time & CirCut & Time \\ & VNSPR & Time & SA & Time\\
    \hline
toursg3-15    &281029888&1023.2&283538879&539.032\\
              &264732800&89044.6&284493491&218.204\\
toursg3-8     &40314704&65.7&41684814&50.527\\
              &41106855&36512.3&41593110&191.472\\
tourspm3-15-50&2964&333.4&2982&411.101\\
              &2834&55083.9&3008&240.492\\
tourspm3-8-50 &442&48.5&454&30.208\\
              &398&9215.2&458&191.124\\
  \end{tabular}
  \caption{Comparison on Set 3 instances.}
  \label{set3}
  \end{figure}

  \begin{figure}
  \begin{tabular}{ll|ll|ll|ll|ll}
    Set & Set size & SS & only & CirCut & only & VNSPR & only & SA & only\\
    \hline
    1 & 54 & 16&9&16&11&9&2&29&19\\
    2 & 30 & 11&0&11&1&5&1&29&16\\
    3 & 4  & 0&0&1&1&0&0&3&3\\
  \end{tabular}
  \caption{For each set, how many times each heuristic discovered the best
    solution among the four.}
  \label{summarytable}
  \end{figure}

  The simulated annealing code in Figure \ref{code} found the best solution
  of the four heuristics on 61 of 89 instances.
  The tables distributed via Optsicom's website appear to list the best
  solution of SS, CirCut, and VNSPR as the best solution known.
  Assuming this is accurate, the code of Figure \ref{code} found solutions
  better than the best previously known solution in 40 of 89 instances.  When
  also comparing against the CirCut runs in this paper, the code of Figure
  \ref{code} found solutions better than any solution found by any other
  heuristic tested in 36 of 89 instances.

  It bears mention that the annealing schedule can be modified in order to get
  much better results on some cases.  For example, changing \texttt{heatmax} to
  40000 and the increment to \texttt{heat} to \texttt{5e-6} finds a solution to
  case g35 from Set 1 of objective value 7685, which is better than those found
  by any other heuristic.  Running with different random seeds can also produce
  better results; inserting the call \texttt{srand(3);} at the beginning of
  \texttt{main} and running on case toursg3-8 from Set 3 gave a solution with
  objective value 41684814, matching CirCut's solution.  I made no special
  effort to optimise the annealing schedule, and the results reported in the
  tables above represent a single run of simulated annealing using the default
  random seed.
\end{section}

\begin{section}{Conclusion}
  The results in this paper appear to show that a simple-minded implementation
  of simulated annealing can often compute within a modest timeframe better
  solutions on widely-known instances than the best algorithms that have been
  applied by the metaheuristics community to the maximum cut problem.  The
  formulation used dates back (at least) to \citet{JohnsonAragon1}.
  This highlights the need to consider simple, classical approaches as well as
  more modern approaches when assessing the state of the art for solving any
  particular problem or establishing a baseline for a computational study.
  
  This should not be considered an endorsement of any particular implementation
  of simulated annealing.  The heuristics involved were run under different
  time limits, with different termination conditions and on different hardware.
  An obvious refinement of this work---for a researcher with access to
  Mart\'{i} et al.'s Scatter Search code---would be to compare all of the
  heuristics considered here on the same computer with the same, fixed, time
  horizon.
\end{section}

\begin{section}{Acknowledgements}
  I thank Levent Tun\c{c}el and Jennifer Yeelam Wong for their comments on
  early versions of this manuscript.  I also thank the Optsicom project for
  collecting and distributing test data and tables in a convenient form.
  This research was supported in part by an NSERC doctoral scholarship.
\end{section}

\renewcommand{\bibname}{References}
\addcontentsline{toc}{chapter}{References}
\bibliographystyle{plainnat}
\bibliography{annealing}
\nocite{*}
\end{document}